\documentclass[12pt]{amsart}
\usepackage{amsmath, amssymb, latexsym, amsthm, mathrsfs}
\usepackage[all]{xy}

      \newenvironment{changemargin}[2]{\begin{list}{}{
         \setlength{\topsep}{0pt}\setlength{\leftmargin}{0pt}
         \setlength{\rightmargin}{0pt}
         \setlength{\listparindent}{\parindent}
         \setlength{\itemindent}{\parindent}
         \setlength{\parsep}{0pt plus 1pt}
         \addtolength{\leftmargin}{#1}\addtolength{\rightmargin}{#2}
         }\item }{\end{list}}

\newcount\skewfactor
\def\mathunderaccent#1#2 {\let\theaccent#1\skewfactor#2
\mathpalette\putaccentunder}
\def\putaccentunder#1#2{\oalign{$#1#2$\crcr\hidewidth
\vbox to.2ex{\hbox{$#1\skew\skewfactor\theaccent{}$}\vss}\hidewidth}}

\newcommand{\arrays}{{{\{0,1\}}^{\N\x\N}}}
\newcommand{\CH}{the Continuum Hypothesis}

\renewcommand{\ss}{\hspace{0.07cm}}
\newcommand{\domi}[3]{[{#1}\ss{#2}\ss{#3}]}
\renewcommand{\mid}[5]{[{#4}\ss{#1}\ss{#3}\ss{#2}\ss{#5}]}

\newcommand{\emp}{excluded middle property}
\newcommand{\scrA}{\mathscr{A}}
\newcommand{\scrB}{\mathscr{B}}

\newcommand{\sr}[2]{{\txt{$#1$\\$#2$}}}
\newcommand{\seq}[1]{{\<#1 : n\in\N\>}}
\newcommand{\setseq}[1]{{\{#1 : n\in\N\}}}

\newcommand{\op}{\operatorname}

\newcommand{\cA}{\mathcal{A}}
\newcommand{\B}{\mathcal{B}}

\newcommand{\BG}{\B_\Gamma}
\newcommand{\BT}{\B_\Tau}
\newcommand{\CG}{C_\Gamma}
\newcommand{\CT}{C_\Tau}
\newcommand{\CO}{C_\Omega}

\newcommand{\BO}{\B_\Omega}

\newcommand{\Tau}{\mathrm{T}}
\newcommand{\cF}{\mathcal{F}}

\newcommand{\M}{\mathcal{M}}
\newcommand{\N}{\mathbb{N}}
\newcommand{\NN}{{\N^\N}}

\renewcommand{\O}{\mathcal{O}}

\newcommand{\cU}{\mathcal{U}}
\newcommand{\Union}{\bigcup}

\newcommand{\Impl}{\Rightarrow}
\long\def\forget#1\forgotten{}

\renewcommand{\b}{\mathfrak{b}}

\renewcommand{\t}{\mathfrak{t}}

\renewcommand{\d}{\mathfrak{d}}

\renewcommand{\i}{\item}

\newcommand{\oo}{\infty}
\newcommand{\p}{\mathfrak{p}}

\newcommand{\fo}{\mathfrak{od}}
\newcommand{\s}{\mathfrak{s}}

\newcommand{\x}{\times}

\newcommand{\nin}{\not\in}

\newcommand{\sbst}{\subseteq}

\newcommand{\sm}{\setminus}

\newcommand{\<}{\langle}
\renewcommand{\>}{\rangle}

\newcommand{\dom}{\op{dom}}
\newcommand{\cov}{\mathsf{cov}}
\newcommand{\add}{\mathsf{add}}

\newcommand{\non}{\mathsf{non}}

\newtheorem{thm}{Theorem}[section]

\newtheorem{lem}[thm]{Lemma}

\newtheorem{cor}[thm]{Corollary}

\theoremstyle{definition}
\newtheorem{defn}[thm]{Definition}
\theoremstyle{remark}

\newcommand{\be}{\begin{enumerate}}
\newcommand{\ee}{\end{enumerate}}
\newcommand{\bi}{\begin{itemize}}
\newcommand{\ei}{\end{itemize}}
\newcommand{\sone}{\mathsf{S}_1}    \newcommand{\sfin}{\mathsf{S}_{fin}}
\newcommand{\ufin}{\mathsf{U}_{fin}}
\newcommand{\newS}[2]{\mathsf{S}_{[#1,#2)}}
\newcommand{\newU}[2]{\mathsf{U}_{[#1,#2)}}

\author{Boaz Tsaban}
\thanks{Partially supported by the Koshland Center for Basic Research.}
\address{Department of Mathematics,
Weizmann Institute of Science, Rehovot 76100, Israel}
\email{boaz.tsaban@weizmann.ac.il}
\urladdr{http://www.cs.biu.ac.il/\~{}tsaban}

\title{Several comments about the combinatorics of $\tau$-covers}

\begin{document}
\begin{abstract}
In a previous work with Mildenberger and Shelah, we showed that the combinatorics of the
selection hypotheses involving $\tau$-covers is sensitive to the
selection operator used.
We introduce a natural generalization of Scheepers' selection operators,
and show that:
\be
\i A slight change in the selection operator, which in classical cases makes no difference,
leads to different properties when $\tau$-covers are involved.
\i One of the newly introduced properties sheds some light on a problem of Scheepers
concerning $\tau$-covers.
\ee
Improving an earlier result, we also show that no generalized Luzin set satisfies $\ufin(\Gamma,\Tau)$.
\end{abstract}


\maketitle

\section{Introduction}

Topological properties defined by diagonalizations of open or Borel covers
have a rich history in various areas of general topology and analysis,
and they are closely related to infinite combinatorial notions,
see \cite{LecceSurvey, futurespm, KocSurv, ict} for surveys on the topic
and some of its applications and open problems.

Let $X$ be an infinite set.
By \emph{a cover of $X$} we mean a family $\cU$ with $X\nin\cU$ and $X=\cup\cU$.
A cover $\cU$ of $X$ is said to be
\be
\i a \emph{large cover of $X$} if: $(\forall x\in X)\ \{U\in\cU : x\in U\}$ is infinite.
\i an \emph{$\omega$-cover of $X$} if: $(\forall\mbox{finite }F\sbst X)(\exists U\in\cU)\ F\sbst U$.
\i a \emph{$\tau$-cover of $X$} if:
$\cU$ is a large cover of $X$, and
$(\forall x,y\in X)\ \{U\in\cU : x\in U\mbox{ and }y\nin U\}$ is finite, or
$\{U\in\cU : y\in U\mbox{ and }x\nin U\}$ is finite.
\i a \emph{$\gamma$-cover of $X$} if: $\cU$ is infinite and
$(\forall x\in X)\ \{U\in \cU : x\nin U\}$ is finite.
\ee
Let $X$ be an infinite, zero-dimensional, separable metrizable topological space
(in other words, a set of reals).
Let $\Omega$, $\Tau$ and $\Gamma$ denote the collections of all \emph{open}
$\omega$-covers, $\tau$-covers and $\gamma$-covers of $X$, respectively.
Additionally, denote the collection of all open covers of $X$ by $\O$.
Similarly, let $\CO$, $\CT$, $\CG$, and $C$ denote the corresponding collections of
\emph{clopen} covers.
Our restrictions on $X$ imply that each member of any of the above classes contains
a countable member of the same class \cite{splittability}.
We therefore confine attention in the sequel to \emph{countable} covers,
and restrict the above four classes to contain only their countable members.
Having this in mind, we let $\BO$, $\BT$, $\BG$, and $\B$ denote the corresponding collections of
\emph{countable Borel} covers.

\begin{figure}[!ht]
\renewcommand{\sr}[2]{{\txt{$#1$\\$#2$}}}
{\tiny
\begin{changemargin}{-3cm}{-3cm}
\begin{center}
$\xymatrix@C=7pt@R=6pt{
%1
&
&
& \sr{\ufin(\Gamma,\Gamma)}{\b}\ar[r]
& \sr{\ufin(\Gamma,\Tau)}{\max\{\b,\s\}}\ar[rr]
&
& \sr{\ufin(\Gamma,\Omega)}{\d}\ar[rrrr]
&
&
&
& \sr{\ufin(\Gamma,\O)}{\d}
\\
%2
&
&
& \sr{\sfin(\Gamma,\Tau)}{\b}\ar[rr]\ar[ur]
&
& \sr{\sfin(\Gamma,\Omega)}{\d}\ar[ur]
\\
%3
\sr{\sone(\Gamma,\Gamma)}{\b}\ar[uurrr]\ar[rr]
&
& \sr{\sone(\Gamma,\Tau)}{\b}\ar[ur]\ar[rr]
&
& \sr{\sone(\Gamma,\Omega)}{\d}\ar[ur]\ar[rr]
&
& \sr{\sone(\Gamma,\O)}{\d}\ar[uurrrr]
\\
%4
&
&
& \sr{\sfin(\Tau,\Tau)}{\min\{\s,\b\}}\ar'[r][rr]\ar'[u][uu]
&
& \sr{\sfin(\Tau,\Omega)}{\d}\ar'[u][uu]
\\
\sr{\sone(\Tau,\Gamma)}{\t}\ar[rr]\ar[uu]
&
& \sr{\sone(\Tau,\Tau)}{\t}\ar[uu]\ar[ur]\ar[rr]
&
& \sr{\sone(\Tau,\Omega)}{\fo}\ar[uu]\ar[ur]\ar[rr]
&
& \sr{\sone(\Tau,\O)}{\fo}\ar[uu]
\\
&
&
& \sr{\sfin(\Omega,\Tau)}{\p}\ar'[u][uu]\ar'[r][rr]
&
& \sr{\sfin(\Omega,\Omega)}{\d}\ar'[u][uu]
\\
\sr{\sone(\Omega,\Gamma)}{\p}\ar[uu]\ar[rr]
&
& \sr{\sone(\Omega,\Tau)}{\p}\ar[uu]\ar[ur]\ar[rr]
&
& \sr{\sone(\Omega,\Omega)}{\cov(\M)}\ar[uu]\ar[ur]\ar[rr]
&
& \sr{\sone(\O,\O)}{\cov(\M)}\ar[uu]
}$
\end{center}
\end{changemargin}
}
\caption{The surviving properties}\label{tauSch}
\end{figure}
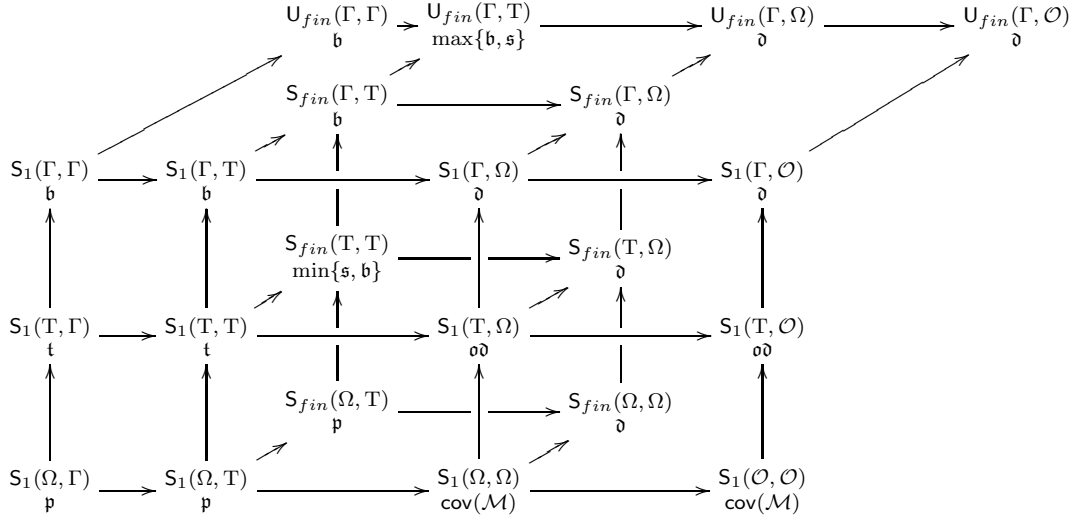

Let $\scrA$ and $\scrB$ be any of the mentioned classes of covers (but of the same descriptive
type, i.e., both open, or both clopen, or both Borel).
Scheepers \cite{coc1} introduced the following \emph{selection hypotheses} that $X$
might satisfy:
\bi
\item[$\sone(\scrA,\scrB)$:]
For each sequence $\seq{\cU_n}$ of members of $\scrA$,
there exist members $U_n\in\cU_n$, $n\in\N$, such that $\setseq{U_n}\in\scrB$.
\item[$\sfin(\scrA,\scrB)$\index{$\sfin(\scrA,\scrB)$}:]
For each sequence $\seq{\cU_n}$
of members of $\scrA$, there exist finite (possibly empty)
subsets $\cF_n\sbst\cU_n$, $n\in\N$, such that $\Union_{n\in\N}\cF_n\in\scrB$.
\item[$\ufin(\scrA,\scrB)$\index{$\ufin(\scrA,\scrB)$}:]
For each sequence $\seq{\cU_n}$ of members of $\scrA$
\emph{which do not contain a finite subcover},
there exist finite (possibly empty) subsets $\cF_n\sbst\cU_n$, $n\in\N$,
such that $\setseq{\cup\cF_n}\in\scrB$.
\ei
Some of the properties are never satisfied, and many equivalences hold among
the meaningful ones. The surviving properties appear in Figure \ref{tauSch},
where an arrow denotes implication \cite{tautau}.
It is not known whether any other implication can be added to this diagram --
see \cite{MShT:858} for a summary of the open problems concerning this diagram.

Below each property $P$ in Figure \ref{tauSch} appears its \emph{critical cardinality},
$\non(P)$, which is the minimal cardinality of a space $X$ not satisfying that property.
The definitions of most of the cardinals appearing in this figure
can be found in \cite{vD, BlassHBK}, whereas $\fo$ is defined in \cite{MShT:858},
and the results were established in
\cite{coc2, tautau, ShTb768, MShT:858}.

A striking observation concerning Figure \ref{tauSch} is, that in the top plane
of the figures, the critical cardinality of $\Pi(\Gamma,\scrB)$ for $\Pi\in\{\sone,\sfin,\ufin\}$
is independent of $\Pi$ in all cases \emph{except for that where $\scrB=\Tau$}.
We demonstrate this anomaly further in Section \ref{gen}, where we also give
a partial answer to a problem of Scheepers.
In Section \ref{luzin} we show that no Luzin set satisfies $\ufin(\Gamma,\Tau)$,
improving a result from \cite{tautau}.

\section{Generalized selection hypotheses}\label{gen}

\begin{defn}
Let $\kappa<\lambda$ be any (finite or infinite) cardinal
numbers. Denote
\bi
\i[$\newS{\kappa}{\lambda}(\scrA,\scrB)$:]
For each sequence $\seq{\cU_n}$ of members of $\scrA$,
there exist subsets $\cF_n\sbst\cU_n$ with $\kappa\le|\cF_n|<\lambda$
for each $n\in\N$, and $\Union_n\cF_n\in\scrB$.
\i[$\newU{\kappa}{\lambda}(\scrA,\scrB)$:]
For each sequence $\seq{\cU_n}$ of members of $\scrA$ which do not contain subcovers
of size less than $\lambda$,
there exist subsets $\cF_n\sbst\cU_n$ with $\kappa\le|\cF_n|<\lambda$
for each $n\in\N$, and $\setseq{\cup\cF_n}\in\scrB$.
\ei
\end{defn}
So that $\newS{1}{2}(\scrA,\scrB)$ is $\sone(\scrA,\scrB)$,
$\newS{0}{\aleph_0}(\scrA,\scrB)$ is $\sfin(\scrA,\scrB)$,
and $\newU{0}{\aleph_0}(\scrA,\scrB)$ is $\ufin(\scrA,\allowbreak\scrB)$.

\begin{defn}\label{semi}
Say that a family $\cA\sbst\arrays$ is \emph{semi $\tau$-diagonalizable} if
there exists a \emph{partial} function $g:\N\to\N$ such that:
\be
\i For each $A\in\cA$: $(\exists^\oo n\in\dom(g))\ A(n,g(n))=1$;
\i For each $A,B\in\cA$:\\
\begin{tabular}{ll}
Either & $(\forall^\oo n\in\dom(g))\ A(n,g(n))\le B(n,g(n))$,\\
or     & $(\forall^\oo n\in\dom(g))\ B(n,g(n))\le A(n,g(n))$.
\end{tabular}
\ee
\end{defn}

In the following theorem, note that $\min\{\s,\b,\fo\}\ge\min\{\s,\b,\cov(\M)\}=\min\{\s,\add\allowbreak(\M)\}$.
\begin{thm}\label{newS}
~\be
\i $X$ satisfies $\newS{0}{2}(\BT,\BT)$ if, and only if, for each Borel function
$\Psi:X\to\arrays$: If $\Psi[X]$ is a $\tau$-family, then it is semi $\tau$-diagonalizable
(Definition \ref{semi}).
The corresponding clopen case also holds.
\i The minimal cardinality of a $\tau$-family that is not semi $\tau$-diagona\-lizable
is at least $\min\{\s,\b,\fo\}$.
\i $\min\{\s,\b,\fo\}\le\non(\newS{0}{2}(\BT,\BT))=\non(\newS{0}{2}(\Tau,\Tau))=\non\allowbreak(\newS{0}{2}(\CT,\CT))$.
\ee
\end{thm}
\begin{proof}
(1) is proved as usual,
(2) is shown in the proof of Theorem 4.15 of \cite{MShT:858},
and (3) follows from (1) and (2).
\end{proof}

\begin{defn}[\cite{ShTb768}]
For functions $f,g,h\in\NN$, and binary relations $R,S$
on $\N$, define subsets $\domi{f}{R}{g}$ and $\mid{R}{S}{g}{h}{f}$
of $\N$ by:
$$\domi{f}{R}{g} = \{ n : f(n)Rg(n)\},\ \mid{R}{S}{g}{f}{h}=\domi{f}{R}{g}\cap\domi{g}{S}{h}.$$
For a subset $Y$ of $\NN$ and $g\in\NN$, we say that
$g$ \emph{avoids middles} in $Y$ with respect to $\<R,S\>$ if:
\be
\i for each $f\in Y$, the set $\domi{f}{R}{g}$ is infinite;
\i for all $f,h\in Y$ at least one of the sets
$\mid{R}{S}{g}{f}{h}$ and $\mid{R}{S}{g}{h}{f}$ is finite.
\ee
$Y$ satisfies the \emph{$\<R,S\>$-\emp{}}
if there exists $g\in\NN$ which avoids middles in $Y$
with respect to $\<R,S\>$.
\end{defn}

In \cite{tautau} it is proved that $\ufin(\BG,\BT)$ is
equivalent to having all Borel images in $\NN$ satisfying
the $\<<,\le\>$-\emp{} (the statement in \cite{tautau} is different
but equivalent).

\begin{thm}\label{combGT}
For a set of reals $X$, the following are equivalent:
\be
\i $X$ satisfies $\newU{1}{\aleph_0}(\BG,\BT)$.
\i Each Borel image of $X$ in $\NN$ satisfies the $\<\le,<\>$-\emp{}.
\ee
The corresponding assertion for $\newU{1}{\aleph_0}(\CG,\CT)$
holds when ``Borel'' is replaced by ``continuous''.
\end{thm}
\begin{proof}
The proof is similar to the one given in \cite{tautau}
for $\ufin(\BG,\BT)$, but is somewhat simpler.

$1\Rightarrow 2$:
Assume that $Y\sbst\NN$ is a Borel image of $X$.
Then $Y$ satisfies $\newU{1}{\aleph_0}(\BG,\BT)$.
For each $n$, the collection $\cU_n=\{U^n_m : m\in\N\}$,
where $U^n_m = \{f\in\NN : f(n)\le m\}$,
is a clopen $\gamma$-cover of $\NN$.
By standard arguments (see $(1\Impl 2)$ in the proof of Theorem 2.3 of \cite{MShT:858})
we may assume that no $\cU_n$ contains a finite cover.
For all $n$, the sequence $\{U^n_m : m\in\N\}$
is monotonically increasing with respect to $\sbst$,
therefore---as large subcovers of $\tau$-covers are also $\tau$-covers---we
may use $\sone(\BG,\BT)$
instead of $\newU{1}{\aleph_0}(\BG,\BT)$ to get a $\tau$-cover
$\cU=\setseq{\Psi^{-1}[U^n_{m_n}]}$
for $X$.
Let $g\in\NN$ be such that $g(n)=m_n$ for all $n$.
Then $g$ avoids middles in $Y$ with respect to $\<\le,<\>$.

$2\Rightarrow 1$:
Assume that $\cU_n=\{U^n_m : m\in\N\}$, $n\in\N$, are Borel covers of $X$
which do not contain a finite subcover.
Replacing each $U^n_m$ with the Borel set $\Union_{k\le m} U^n_k$ we may assume
that the sets $U^n_m$ are monotonically increasing with $m$.
Define $\Psi:X\to\NN$ by: $\Psi(x)(n) = \min\{m : x\in U^n_m\}$.
Then $\Psi$ is a Borel map, and so $\Psi[X]$
satisfies the $\<\le,<\>$-\emp{}.
Let $g\in\NN$ be a witness for that.
Then $\cU=\setseq{U^n_{g(n)}}$ is a $\tau$-cover of $X$.

The proof in the clopen case is similar.
\end{proof}

\begin{cor}\label{anomaly}
The critical cardinalities of $\newU{1}{\aleph_0}(\BG,\BT)$, $\newU{1}{\aleph_0}(\Gamma,\Tau)$, and
$\newU{1}{\aleph_0}(\CG,\allowbreak\CT)$, are all equal to $\b$.
\end{cor}
\begin{proof}
This follows from Theorem \ref{combGT} and the corresponding combinatorial assertion,
which was proved in \cite{ShTb768}.
\end{proof}

Recall from Figure \ref{tauSch} that the critical cardinality of
$\ufin(\Gamma,\Tau)=\newU{0}{\aleph_0}(\Gamma,\Tau)$ is $\max\{\s,\b\}$.
Contrast this with Corollary \ref{anomaly}.

According to Scheepers \cite[Problem 9.5]{futurespm},
one of the more interesting problems concerning Figure \ref{tauSch} is
whether $\sone(\Omega,\Tau)$ implies $\ufin(\Gamma,\Gamma)$.
If $\newU{1}{\aleph_0}(\Gamma,\Tau)$ is preserved under taking finite unions,
then we get a positive solution to Scheepers' Problem.
(Note that $\sone(\Omega,\Tau)$ implies $\sone(\Gamma,\Tau)$.)

\begin{cor}
If $\newU{1}{\aleph_0}(\Gamma,\Tau)$ is preserved under taking finite unions,
then it is equivalent to $\ufin(\Gamma,\Gamma)$ and
$\sone(\Gamma,\Tau)$ implies $\ufin(\Gamma,\Gamma)$.
\end{cor}
\begin{proof}
The last assertion of the theorem follows from the first since
$\sone(\Gamma,\Tau)$ implies $\newU{1}{\aleph_0}(\Gamma,\Tau)$.

Assume that $X$ does not satisfy $\ufin(\Gamma,\Gamma)$.
Then, by Hurewicz' Theorem \cite{HURE27},
there exists an unbounded continuous image $Y$ of $X$ in $\NN$.
For each $f\in Y$, define $f_0,f_1\in\NN$ by $f_i(2n+i)=f(n)$ and $f_i(2n+(1-i))=0$.
For each $i\in\{0,1\}$, $Y_i = \{f_i : f\in Y\}$ is a continuous image of $Y$.
It is not difficult to see that $Y_0\cup Y_1$ does not satisfy the  $\<\le,<\>$-\emp{}
\cite{ShTb768}.
By Theorem \ref{combGT}, $Y_0\cup Y_1$ does not satisfy $\newU{1}{\aleph_0}(\Gamma,\Tau)$,
thus, by the theorem's hypothesis, one of the sets $Y_i$ does not satisfy that property.
Therefore $Y$ (and therefore $X$) does not satisfy $\newU{1}{\aleph_0}(\Gamma,\Tau)$ either.
\end{proof}

We do not know whether $\newU{1}{\aleph_0}(\Gamma,\Tau)$ is
preserved under taking finite unions.
We also do not know the situation for $\ufin(\Gamma,\Tau)$.
The following theorem is only interesting when $\s<\b$.

\begin{thm}
If there exists a set of reals $X$ satisfying $\ufin(\Gamma,\Tau)$ but not $\ufin(\Gamma,\Gamma)$,
then $\ufin(\Gamma,\Tau)$ is not preserved under taking unions of $\s$ many elements.
\end{thm}
\begin{proof}
The proof is similar to the last one, except that here we define $\s$ many continuous images
of $Y$ as we did in \cite{ShTb768} to prove that the critical cardinality of
$\ufin(\Gamma,\Tau)$ is $\max\{\s,\b\}$.
\end{proof}

\section{Luzin sets}\label{luzin}

A set of reals $L$ is a \emph{generalized Luzin set} if for each meager set $M$,
$|L\cap M|<|L|$.
In \cite{tautau} we constructed (assuming a portion of \CH{})
a generalized Luzin set which satisfies
$\sone(\BO,\BO)$ but not $\ufin(\Gamma,\Tau)$.
We now show that the last assertion always holds.

\begin{thm}
Assume that $L\sbst\NN$ is a generalized Luzin set.
Then $L$ does not satisfy the $\<<,\le\>$-\emp{}.
In particular, $L$ does not satisfy $\ufin(\CG,\CT)$.
\end{thm}
\begin{proof}
We use the following easy observation.
\begin{lem}[\cite{tautau}]
Assume that $A$ is an infinite set of natural numbers, and $f\in\NN$.
Then the sets
\begin{eqnarray*}
       M_{f,A} & = & \{ g\in\NN : [g\le f]\cap A\mbox{ is finite}\}\\
\tilde M_{f,A} & = & \{ g\in\NN : [f < g]\cap A\mbox{ is finite}\}
\end{eqnarray*}
are meager subsets of $\NN$.\hfill\qed
\end{lem}
Fix any $f\in\NN$. We will show that $f$ does not avoid middles in $Y$ with
respect to $\<<,\le\>$.
The sets
$M_{f,\N} = \{ g\in\NN : [g\le f]\mbox{ is finite}\}$ and
$\tilde M_{f,\N} = \{ g\in\NN : [f < g]\mbox{ is finite}\}$
are meager, thus there exists $g_0\in L\sm (M_{f,\N}\cup\tilde M_{f,\N})$.
Now consider the meager sets
$M_{f,[f < g_0]} = \{ g\in\NN : [g\le f < g_0]\mbox{ is finite}\}$ and
$\tilde M_{f,[g_0\le f]} = \{ g\in\NN : [g_0\le f<g]\mbox{ is finite}\}$,
and choose $g_1\in L\sm(M_{f,[f < g_0]}\cup \tilde M_{f,[g_0\le f]})$.
Then both sets $[g_0<f\le g_1]$ and $[g_1<f\le g_0]$ are infinite.
\end{proof}

\end{document}